\documentclass[11pt]{article}
\usepackage{graphicx}
\usepackage{mathrsfs}
\usepackage{latexsym}
\usepackage{psfrag}
\usepackage{graphicx,subfigure}
\usepackage{fancyhdr,graphicx}
\usepackage{multicol}
\usepackage{cite}
\usepackage{amsmath,amssymb}
\usepackage{flafter}
\usepackage{fancyhdr}
\usepackage{subfigure}
\usepackage{indentfirst,latexsym,bm}
\usepackage{xcolor}
\usepackage{float}
\usepackage{appendix}

\usepackage{tcolorbox}
\usepackage{hyperref}
\hypersetup{colorlinks=true,linkcolor=blue,filecolor=blue,citecolor=blue}
\usepackage{booktabs}

\usepackage{algorithm}  
\usepackage{algpseudocode}  
\usepackage{amsmath}  

\newtheorem{theorem}{Theorem}[section]
\newtheorem{lemma}[theorem]{Lemma}

\newtheorem{conjecture}[theorem]{Conjecture}

\newenvironment{pf}{
\par
\noindent {\bf Proof.}\rm}
{\mbox{}\hfill\rule{0.5em}{0.8em} \par \bigskip}

\begin{document}

\title{\bf The tight bound for the strong chromatic indices of claw-free subcubic graphs
\thanks{Supported by NSFC 11771080.}}
\author{Yuquan Lin and Wensong Lin\footnote{Corresponding author. E-mail address: wslin@seu.edu.cn}\\
{\small School of Mathematics, Southeast University, Nanjing 210096, P.R. China}}
\date{}
\maketitle

\vspace*{-1cm} \setlength\baselineskip{6.5mm}
\bigskip

\begin{abstract}

Let $G$ be a graph and $k$ a positive integer. A strong $k$-edge-coloring of $G$ is a mapping $\phi: E(G)\to \{1,2,\dots,k\}$ such that for any two edges $e$ and $e'$ that are either adjacent to each other or adjacent to a common edge, $\phi(e)\neq \phi(e')$. The strong chromatic index of $G$,  denoted as $\chi'_{s}(G)$,  is the minimum integer $k$ such that $G$ has a strong $k$-edge-coloring. Lv, Li and Zhang [Graphs and Combinatorics 38 (3) (2022) 63] proved that if $G$ is a claw-free subcubic graph other than the triangular prism then $\chi_s'(G)\le 8$. In addition, they asked if the upper bound $8$ can be improved to $7$. In this paper, we answer this question in the affirmative. Our proof implies a linear-time algorithm for finding strong $7$-edge-colorings of such graphs. We also construct infinitely many claw-free subcubic graphs with their strong chromatic indices attaining the bound $7$.\\
 

\end{abstract}

\section{Introduction}

Let $G=(V(G), E(G))$ be a finite undirected simple  graph.
For $v\in V(G)$, let $N(v)=\{u\in V(G): uv\in E(G)\}$ denote the open neighborhood of $v$ and  $d(v)=|N(v)|$ be the degree of $v$. Let  $\Delta(G)=\max\limits_{v\in V(G)}d(v)$ denote the maximum degree of $G$.  
For convenience, we use the abbreviation $[1,n]$ for $\{1,2, \dots, n\}$, where $n$ is any positive integer. 

Let $e$ and $e'$ be two edges of $G$. If $e$ and $e'$ are adjacent to each other, we say that the distance between $e$ and $e'$ is $1$, and if they are not adjacent but both of them are adjacent to a common edge, we say they are at distance $2$. Given a positive integer $k$,  a {\em strong $k$-edge-coloring} of $G$ is a mapping $\phi: E(G)\to [1,k]$ such that for any two edges $e$ and $e'$ that are at distance $1$ or $2$, $\phi(e)\neq \phi(e')$. The {\em strong chromatic index} of $G$,  denoted by $\chi'_{s}(G)$,  is the minimum integer $k$ such that $G$ has a strong $k$-edge-coloring.

The concept of strong edge coloring, first introduced by Fouquet and Jolivet \cite{FJ1983}, can be used to model the conflict-free channel assignment in radio networks \cite{R1997,NKGB2000}. In 1985,   Erd\H{o}s and  Ne\v{s}et\v{r}il \cite{E1988,EN1989} proposed the following conjecture about the upper bound of $\chi'_{s}(G)$ in term of $\Delta(G)$, which if true, is the best possible.

\begin{conjecture}\label{Conj-EN}
	(Erd\H{o}s and Ne\v{s}et\v{r}il \cite{E1988,EN1989})
	If $G$ is a graph with maximum degree $\Delta(G)$, then
	\begin{equation*}
	\chi'_{s}(G)\le\begin{cases}
			\begin{array}{cl}
		\dfrac{5}{4}\Delta(G)^{2},& \text{if}\  \Delta(G) \ \text{is even,} \\
			\dfrac{5}{4}\Delta(G)^{2}-\dfrac{1}{2}\Delta(G)+\dfrac{1}{4},& \text{if}\  \Delta(G)\  \text{is odd.}
				\end{array}
	\end{cases}	
	\end{equation*}
\end{conjecture}

The conjecture is clearly true for $\Delta(G)\le 2$. The case $\Delta(G)=3$ was verified  by Andersen \cite{A1992} in 1992, and independently by Hor\'{a}k, Qing, and Trotter \cite{HQT1993} in 1993. Furthermore, if $G$ is a subcubic planar graph then $\chi'_s(G)\le 9$ \cite{KLRSWY2016}. For  $\Delta(G)=4$,  an upper bound of $23$ was proved by Hor\'{a}k \cite{H1990} in 1990. It was improved to $22$ by  Cranston \cite{C2006} in 2006 and more recently Huang,  Santana and Yu \cite{HSY2018} obtained the upper bound $21$.  Recall that the conjectured bound is $20$ for $\Delta(G)=4$. Although the general case when $\Delta(G)=4$ has not been solved yet, the bound $20$ holds for some special cases. In 1990, Faudree, Schelp, Gy{\'a}rf{\'a}s  and Tuza \cite{FSGT1990} proved that $\chi'_{s}(G)\le  4\Delta(G) + 4$ for any planar graph $G$, which implies that Conjecture \ref{Conj-EN} holds for planar graphs with maximum degree $4$. And this bound for planar graphs with maximum degree $4$ was further improved to $19$ by Wang, Shiu, Wang and Chen \cite{WSWC2018} in 2018.

In 2015, Bensmail, Bonamy and  Hocquard \cite{BBH2015} studied the upper bound on $\chi'_{s}(G)$  in terms of maximum average degree for the class of  graphs with maximum degree $4$, where they gave some  sufficient conditions under which Conjecture \ref{Conj-EN} is true. In 2018, Lv, Li and Yu \cite{LLY2018} strengthened Bensmail, Bonamy and  Hocquard's results. And they  proved that for any  graph $G$ with $\Delta(G)=4$, if  there are two vertices of degree $3$ whose distance is at most $4$, then  $\chi'_{s}(G)\le  20$.

For graphs with maximum degree $5$,  Zang \cite{Z2015} gave the upper bound 37 (recall that the conjectured bound is 29), which is the only progress as we know. And for larger $\Delta(G)$, the problem is widely open. In 1997, Molloy and Reed \cite{MR1997} used probabilistic techniques to  prove that $\chi'_{s}(G)\le 1.998\Delta(G)^{2}$ for sufficiently large $\Delta(G)$. And an improvement to $1.93\Delta(G)^{2}$ was provided by Bruhn and Joos \cite{BJ2015}  in 2015. Recently, Bonamy, Perrett and Postle \cite{BPP2022} further strengthened this bound to $1.835\Delta(G)^{2}$. The current best known upper bound is $1.772\Delta(G)^{2}$ which was showed by Hurley, de Joannis de Verclos and Kang \cite{HdK2021} in 2021.

A graph is called {\em claw-free} if it has no induced subgraph isomorphic to the complete bipartite graph $K_{1,3}$. In 2020, Debski, Junosza-Szaniawski and {\'S}leszy{\'n}ska-Nowak  \cite{DJS2020} presented the following upper bound for the strong chromatic indices of claw-free graphs.

\begin{theorem}\label{claw-free}
(D{e}bski, Junosza-Szaniawski and {\'S}leszy{\'n}ska-Nowak \cite{DJS2020})
For any claw-free graph $G$ with maximum degree $\Delta(G)$, $\chi'_{s}(G)\le  \frac{9}{8}\Delta(G)^{2}+\Delta(G)$.
\end{theorem}	

\begin{figure}[htbp]  
	\centering
	\begin{minipage}[t]{5cm}
		\centering
		\resizebox{4.5cm}{3cm}{\includegraphics{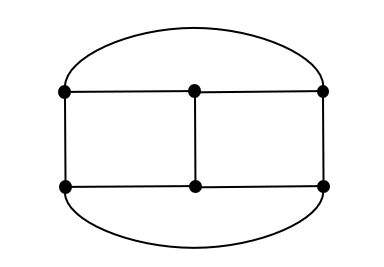}}
		\caption{The $3$-prism}
		\label{fig:3-prism}
	\end{minipage}
	\begin{minipage}[t]{5cm}
		\centering
		\resizebox{3.6cm}{3cm}{\includegraphics{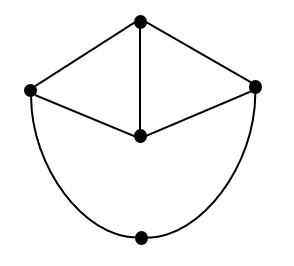}}
		\caption{The graph $G_0$}
		\label{fig:G0}
	\end{minipage}
	\begin{minipage}[t]{5cm}
		\centering
		\resizebox{2.8cm}{3cm}{\includegraphics{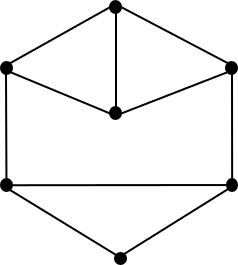}}
		\caption{The graph $H_{0}$}
		\label{fig:H0}
	\end{minipage}
\end{figure}

A graph with maximum degree less than or equal to $3$ is  called a {\em subcubic} graph. In 2022,  Lv, Li and Zhang \cite{LLZ2022} proved that, for any claw-free subcubic graph $G$ other than the triangular prism,  $\chi'_{s}(G)\le  8$. Please see Figure \ref{fig:3-prism} for the triangular prism (also called the $3$-prism).  Notice that the $3$-prism is a claw-free cubic graph with its strong chromatic index being equal to $9$. In the same paper, the authors left the problem whether this bound can be improved to $7$. This paper solves this problem and the main result is the following theorem.

\begin{theorem}\label{main}
Let $G$ be a claw-free subcubic graph. If each component of $G$ is not isomorphic to the triangular prism, then $\chi'_{s}(G)\leq 7$. 	
\end{theorem}

\noindent{\bf Remark 1.}
Recall that when $\Delta (G)= 3$, the upper bound in Conjecture \ref{Conj-EN} is $10$, the
upper bound in Theorem \ref{claw-free} is $\frac{105}{8}$, and the upper bound proved by Lv, Li and Zhang \cite{LLZ2022} is $8$, while our bound in Theorem \ref{main} is $7$. \\

\noindent{\bf Remark 2.} A graph $G_0$ on five vertices with strong chromatic index $7$ was presented in \cite{LLZ2022}  (see Figure \ref{fig:G0}). This shows the sharpness of the upper bound $7$. In fact, there are other claw-free subcubic graphs with their strong chromatic indices being equal to $7$. It is not difficult to verify that the graph $H_{0}$ shown in Figure \ref{fig:H0} has the strong chromatic index $7$. Therefore, by Theorem \ref{main}, any claw-free subcubic graph containing $H_0$ as a subgraph has the strong chromatic index $7$. This implies that there are infinitely many claw-free subcubic graphs with their strong chromatic indices attaining the upper bound $7$.

\noindent{\bf Remark 3.} We would like to point out that our proof of Theorem \ref{main} implies a linear-time algorithm that can produces a strong $7$-edge-coloring of any claw-free subcubic graph other than the $3$-prism.\\

The remainder of this paper is dedicated to the proof of Theorem \ref{main} which is organized as follows. In Section 2, after stating some definitions and notations, we explain how to order the edges of a connected claw-free subcubic graph and apply the greedy algorithm to obtain a partial strong $7$-edge-coloring of the graph. Section 3 consists of a series of lemmas, each of which shows that a particular partial strong $7$-edge-coloring of a claw-free subcubic graph can be extended to a strong $7$-edge-coloring of the whole graph. Finally, Section 4 summarizes our works and suggests some future research directions.

\section{Preliminaries and notations}

Let $G$ be a connected claw-free subcubic graph not isomorphic to the $3$-prism. To prove Theorem \ref{main}, we shall first find a subgraph of $G$ and construct a strong $7$-edge-coloring of this subgraph, and then try to extend the coloring to the whole graph $G$. The subgraph we find should have the property that any two of its edges at distance greater than $2$ must also be at distance greater than $2$ in $G$. We next introduce some notations and preliminary facts that we will use in our proofs.

We use $\alpha$, $\beta$, $\gamma$ to denote colors and $\phi$, $\psi$, $\sigma$ to denote edge colorings. Given two distinct edges $e$ and $e'$ of $G$, we say that $e$ {\em sees} $e'$ in $G$ if they are distance $1$ or $2$ apart. An edge coloring of a graph $G$ is {\em good}, if it is a strong edge coloring of $G$ using at most $7$ colors. A {\em good partial coloring} of a graph $G$ is a good coloring $\phi$ of some subgraph $H$ of $G$ such that $\phi(e)\neq \phi(e')$ if $e$ and $e'$ see each other in $G$.

Let $\phi$ be  a good partial coloring of $G$. We say that $e$ {\em sees} a color $\alpha$ in $\phi$, if $e$ sees an edge $e'$ for which $\phi(e) = \alpha$. For $e\in E(G)$, let $F_{\phi}(e)$ denote the set of colors that  $e$ sees in  $\phi$. In addition, let $\bar{E}_{\phi}$ denote the  set of  edges in $G$ not already assigned colors by $\phi$ and $\bar{G}_{\phi}$ denote the subgraph of $G$ induced by   $\bar{E}_{\phi}$. For  $e\in \bar{E}_{\phi}$, let $A_{\phi}(e)$ denote the set of colors that $e$ does not see in $\phi$.  It is clear that $A_{\phi}(e)=[1,7]\setminus F_{\phi}(e)$ for any $e\in \bar{E}_{\phi}$.

In order to construct a good partial  coloring of a claw-free subcubic graph $G$, we first apply the similar method in \cite{WL2008} to order the edges of $G$, and then use {\em the greedy algorithm} to color them in the order, left only a few particular edges uncolored.   

Suppose $S$ is a nonempty subset of $V(G)$. For a vertex $v\in V(G)$, the distance from $v$ to $S$, denoted by $d_S(v)$, is equal to $\min\limits_{w\in S}\{d(v,w)\}$, where $d(v,w)$ is the distance between $v$ and $w$ in $G$. Let $I$ be the maximum distance from a vertex of $G$ to $S$. For $i=0,1,\dots,I$, let $D_i=\{v\in V(G):\ d(v,S)=i\}$. A mapping $d_S$ from $E(G)$ to nonnegative real numbers is defined as: for any edge $e$ with two end vertices $u$ and $v$, $d_S(e)=\frac{1}{2} (d_S(u)+ d_S(v))$. Suppose $R=(e_{k_1},e_{k_2},\dots,e_{k_m})$ is an ordering of the edges of $G$. For any two integers $i$ and $j$ in $[1,m]$, if $i<j$ implies $d_S(e_{k_i})\geq d_S(e_{k_j})$, then we say that the edge ordering $R$ of $G$ is {\em compatible} with the mapping $d_S$. It is clear that any edge joins two vertices that are either in the same $D_{i}$, or one in $D_{i}$ and the other in $D_{i+1}$ for some $i$, and that each vertex in  $D_{i}$ with $i \ge 1$ is adjacent to at least one vertex in $D_{i-1}$.  Let $e=xy$ be an edge with $d_S(x)\le d_S(y)$. If $d_S(e)\ge 1$, then $e$ is adjacent to an edge $e'=xz$ with $d_S(z)=d_S(x)-1$.

We are now ready to produce a good partial coloring of a claw-free subcubic graph.

\begin{lemma}\label{greedy}
Let $G$ be a connected claw-free subcubic graph other than the $3$-prism and $S$ a nonempty subset of $V(G)$. The greedy algorithm, coloring the edges of $G$ in an order $R$ that is compatible with the  mapping $d_S$, will produce a good partial  coloring of $G$ with only edges in $\{e\in E(G): d_S(e)<1\}$ being left uncolored.
\end{lemma}

\begin{pf}
Let $e=xy$ be an edge with $d_S(e)\ge 1 $. Suppose $d_S(x) \le d_S(y)$. Let $e'=xz$ be a neighbor of $e$ with $d_S(z) = d_S(x)-1$. It is clear that all edges incident to the vertex $z$ are behind the edge $e$ in the order $R$. And so, when $e$ is going to be colored by the greedy algorithm, all edges incident to $z$ are not colored. 

If $d(x)=2$ or $d(x)=3$ and the third neighbor $w$ of $x$ satisfies $d_S(w)=d_S(x)-1$, then $e$ sees at most $5$ colors at that moment, implying that $e$ can be colored properly. Thus we assume that $d(x)=3$ and the third neighbor $w$ of $x$ satisfies $d_S(w) \ge d_S(x)$. 

If $d_S(y)=d_S(x)\ge 1$, then $y$ has a neighbor $y'$ with $d_S(y')=d_S(z)$ (It is possible that $z=y'$). As $G$ is a claw-free subcubic graph, it is easy to check that $e$ sees at most $6$ colors (please refer to  Figure \ref{fig:greedy} (a) and (b)). Thus  $e$ can be colored properly.
 	
		\begin{figure}[htbp]
		\centering
		\resizebox{12cm}{4.6cm}{\includegraphics{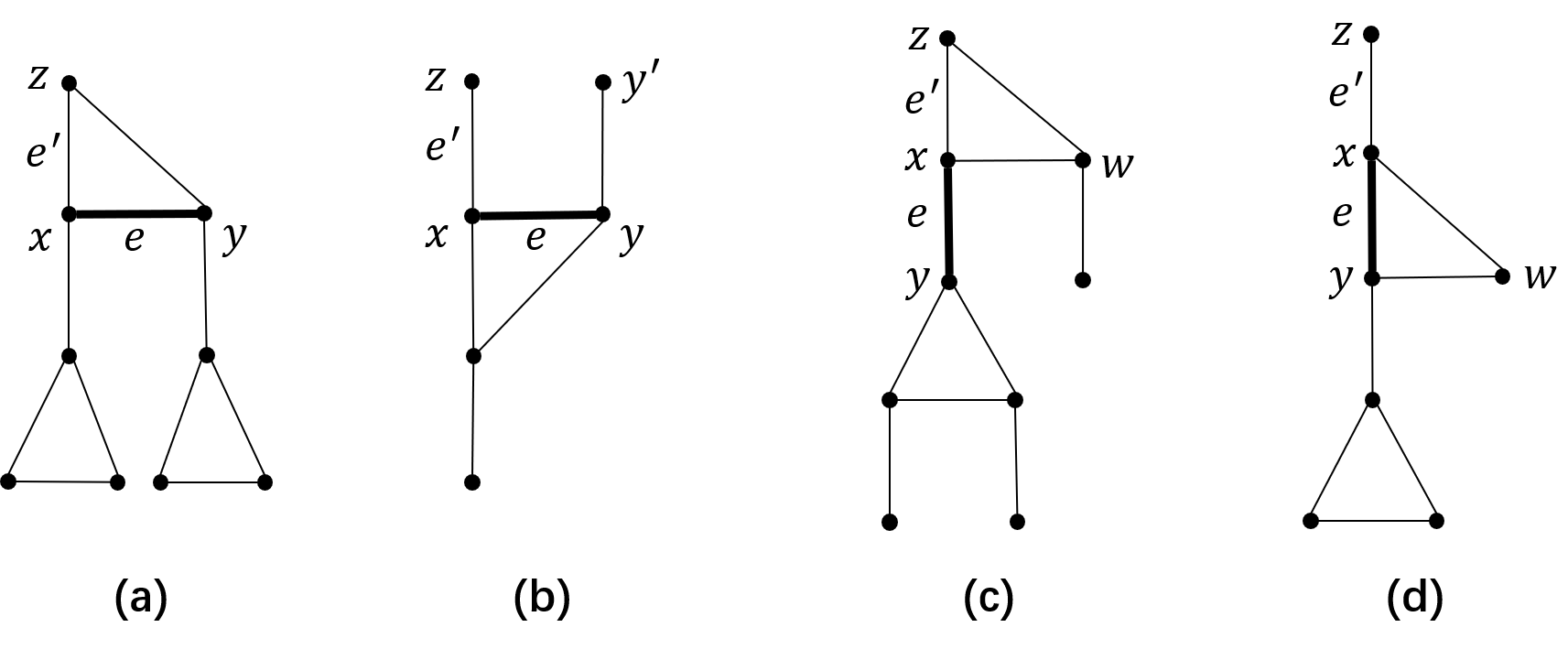}}
		\caption{The illustrations of Lemma \ref{greedy}}
		\label{fig:greedy}
	       \end{figure}
	       	
If $d_S(y)=d_S(x)+1\ge 2$, then $d_S(w)=d_S(x)$ or $d_S(w)=d_S(x)+1$. Again, because $G$ is a claw-free subcubic graph,  $e$ sees at most $6$ colors (please refer to  Figure \ref{fig:greedy} (c) and (d))  and so  $e$ can be colored properly. Therefore, the lemma holds.
\end{pf}

In the proof of Theorem \ref{main}, we will often choose some subset $S$ of $V(G)$. With this subset $S$, by Lemma \ref{greedy}, we can obtain a good partial coloring $\phi$ of $G$ with only edges $e$ with $d_S(e)<1$ being uncolored. And then our task is to extend $\phi$ to a good coloring of the whole graph $G$.

There are two main techniques in our proofs, one is using Hall's theorem \cite{H1935}, the other is the  modified caterpillar tree method which is slightly different from the  caterpillar tree method used by Hor\'{a}k, Qing and Trotter \cite{HQT1993} in proving  the upper bound $10$ for the strong chromatic indices of cubic graphs.

By Hall's theorem \cite{H1935}, $\{A_{\phi}(e):e\in\bar{E}_{\phi}\}$ has a system of distinct representatives  (abbreviated SDR) if and only if $|\cup_{e\in M}A_{\phi}(e)|\ge |M|$, for every $M\subseteq \bar{E}_{\phi}$. Whenever  $\{A_{\phi}(e):e\in\bar{E}_{\phi}\}$ has a SDR, $\phi$ can be easily extended to a good coloring  of $G$. In this situation, we will say that we can obtain a good coloring of $G$ by SDR.

\section{The proof of Theorem  \ref{main} }

It is sufficient to prove Theorem \ref{main} for connected graphs. 
Let $G$ be a connected claw-free subcubic graph not isomorphic to the $3$-prism.
If $G$ is isomorphic to the complete graph $K_{4}$ or the graph $K_{4}^{\Delta}$ (i.e., the graph obtained from $K_{4}$ by replacing each vertex with a $3$-cycle, as shown in Figure \ref{fig:K_{4}}), then it is easy to check that  $\chi'_{s}(G)\le 7$. Please refer to  Figure \ref{fig:K_{4}} for a good coloring of $K_{4}^{\Delta}$. Thus in this section,  we always assume that $G$ is a connected claw-free subcubic graph that is not isomorphic to any graph in $\{\mbox{the $3$-prism}, K_{4}, K_{4}^{\Delta}\}$.

\begin{figure}[htbp]
\centering
\resizebox{14cm}{3.7cm}{\includegraphics{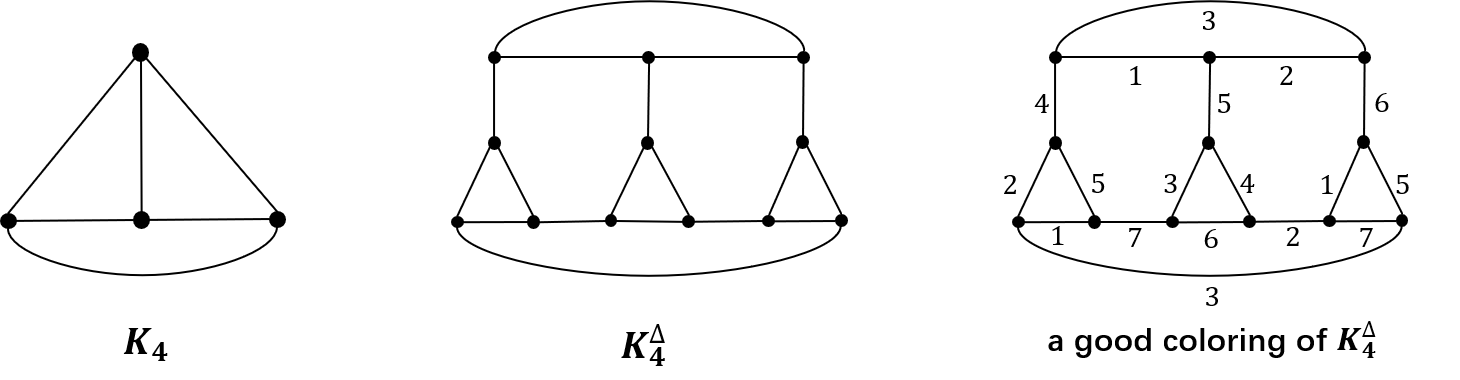}}
\caption{The two graphs $K_{4},K_{4}^{\Delta}$ and a good coloring of $K_{4}^{\Delta}$}
\label{fig:K_{4}}
\end{figure}
    
\begin{lemma}\label{lemma:degree-1}
 	If $G$ has a vertex of degree $1$, then $\chi'_{s}(G)\leq 7$.
\end{lemma}

\begin{pf}	
Let $v_0$ be a vertex of degree $1$ in $G$ and  $e_0$ be the edge incident with $v_0$. 
Set $S=\{v_0\}$. Then,  by Lemma  \ref{greedy}, $G-v_{0}$ has a good coloring.
Since $e_{0}$ sees at most $5$ edges in $G$,  we can extend this good coloring of $G-v_{0}$ to $G$.
\end{pf}
 
\begin{lemma}\label{lemma:degree-2}
 	If $G$ has a vertex of degree $2$, then $\chi'_{s}(G)\leq 7$.
\end{lemma}

\begin{pf} 
Let $v_0$ be a vertex of degree $2$ in $G$ with two neighbors $v_{1}$ and $v_{2}$. If $v_{1}v_{2}\in E(G)$ or $d(v_{1})=2$, then put  $S=\{v_0\}$. By Lemma \ref{greedy}, $G-v_{0}$ has a good coloring $\phi$ with $|A_{\phi}(v_{0}v_{1})|\ge 2$ and $|A_{\phi}(v_{0}v_{2})|\ge 1$. We can obtain a good coloring of $G$ by SDR easily. So by  symmetry, we may assume that $v_{1}v_{2}\notin E(G)$ and $d(v_{1})=d(v_{2})=3$. Let $N(v_{1})=\{v_{0},u_{1},u_{1}'\}$ and $N(v_{2})=\{v_{0},u_{2},u_{2}'\}$. Since $G$ is  a claw-free subcubic graph, we must have $u_1u_1'\in E(G)$ and $u_2u_2'\in E(G)$. This implies that $\{u_1,u_1'\}=\{u_2,u_2'\}$ or $\{u_1,u_1'\}\cap \{u_2,u_2'\}=\emptyset$.

If $\{u_1,u_1'\}=\{u_2,u_2'\}$ then $G$ is isomorphic to $H_{1}$ (please see Figure \ref{fig:H_{1}H_{2}H_{3}} for $H_{1}$), and so $\chi'_{s}(G)= 7$. Thus we assume that $\{u_1,u_1'\}\cap \{u_2,u_2'\}=\emptyset$.  If $d(u_{1})=2$, then the two neighbors of $u_1$ are $v_1$ and $u_1'$. By setting $S=\{u_1\}$, similar to the argument in the previous paragraph, we can get a good coloring of $G$. Thus we assume $d(u_{1})=3$. Symmetrically, we may also assume that $d(u_{1}')=d(u_{2})=d(u_{2}')=3$. Let $w_{1},w_{1}',w_{2},w_{2}'$ denote the third neighbors of $u_{1},u_{1}',u_{2},u_{2}'$, respectively. Denote by $E(\{u_1,u_1'\},\{u_2,u_2'\})$ the set of edges with one vertex in $\{u_1,u_1'\}$ and the other in $\{u_2,u_2'\}$. If $|E(\{u_1,u_1'\},\{u_2,u_2'\})|=2$, then $G$ is isomorphic to $H_2$ (see Figure \ref{fig:H_{1}H_{2}H_{3}}) and we clearly have $\chi'_{s}(G)\le 7$. We next deal with the remaining two cases $|E(\{u_1,u_1'\},\{u_2,u_2'\})|=1$ and $|E(\{u_1,u_1'\},\{u_2,u_2'\})|=0$.
	
 		\begin{figure}[htbp]
 		\centering
 		\resizebox{9.5cm}{2.7cm}{\includegraphics{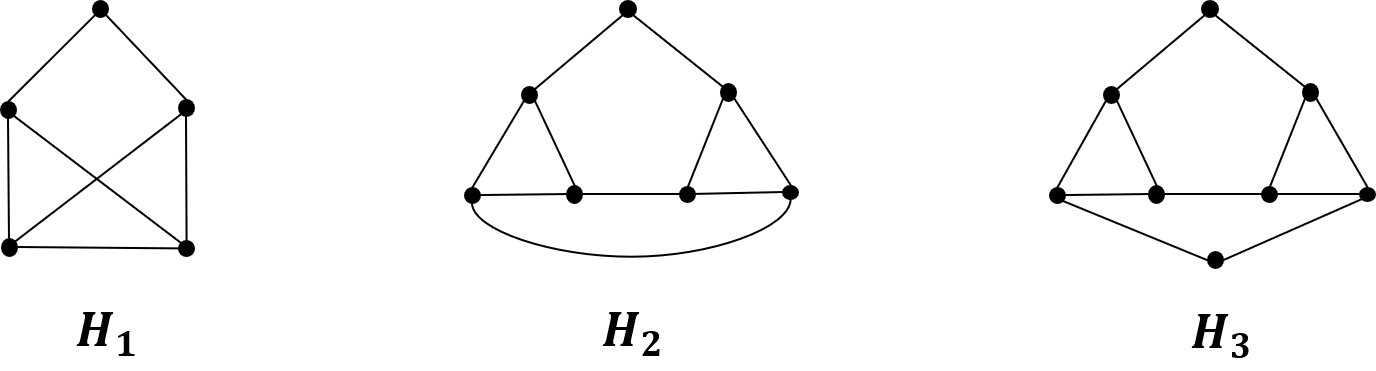}}
 		\caption{The three graphs $H_{1},H_{2}$ and $H_{3}$}
 		\label{fig:H_{1}H_{2}H_{3}}
 	     \end{figure} 
	     
{\bf Case 1}. $|E(\{u_1,u_1'\},\{u_2,u_2'\})|=1$.

Without loss of generality, assume that $u_{1}u_{2}\in E(G)$. If $u_{1}'$ and $u_{2}'$ are adjacent to a common neighbor, then $G$ must be isomorphic to $H_{3}$ (see Figure \ref{fig:H_{1}H_{2}H_{3}}), implying that $\chi'_{s}(G)\le 7$.  Thus, we assume $u_{1}'u_{2}'\notin E(G)$ and $w_{1}'\neq w_{2}'$.

Set $S=\{v_{0},v_{1},v_{2},u_{1},u_{2}\}$. By  Lemma  \ref{greedy}, we get a good partial coloring $\phi$ of $G$ with nine edges uncolored (please see Figure \ref{fig:degree2_case1} for the names of the uncolored edges of $G$). Observe that $|A_{\phi}(e_{3})|\ge5$, $|A_{\phi}(e_{i})|=6$ for $i=1,2,4,5$ and $|A_{\phi}(f_{j})|\ge 4$ for $j=1,2,3,4$.

	\begin{figure}[htbp]
	\centering
	\resizebox{11.1cm}{4cm}{\includegraphics{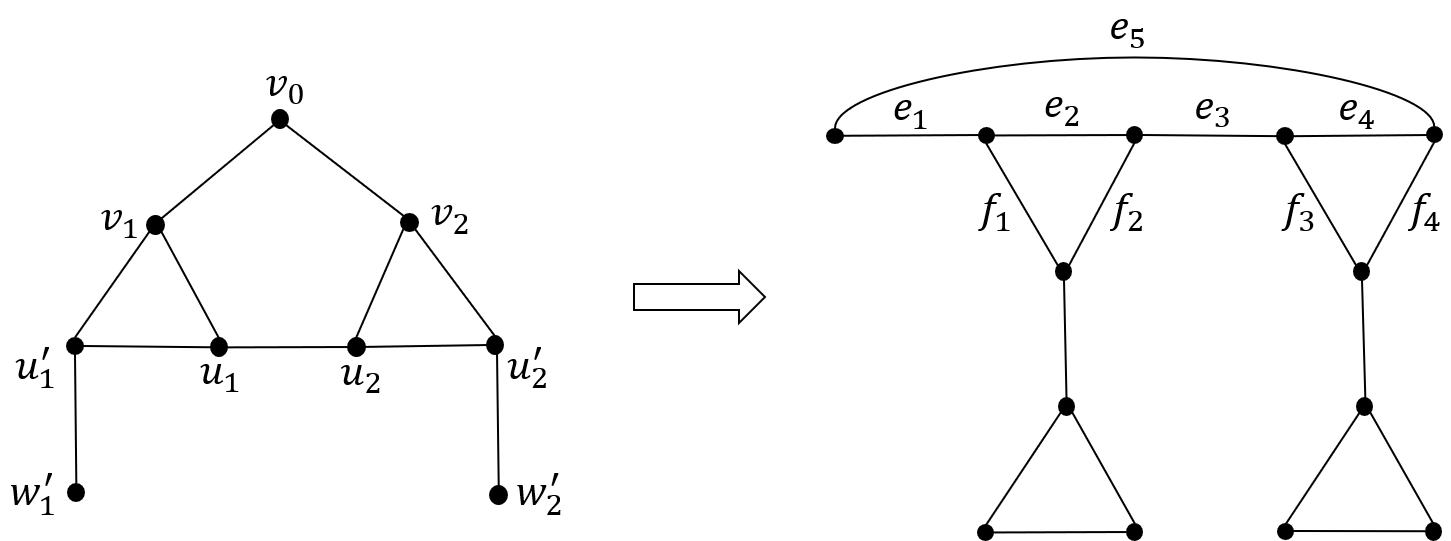}}
	\caption{Case 1 in the proof of Lemma \ref{lemma:degree-2}}
	\label{fig:degree2_case1}
       \end{figure}

Let $\alpha$ be a color in $A_{\phi}(f_{1})\cap A_{\phi}(f_{3})$. By coloring $f_{1}$ and $f_{3}$ with the same color $\alpha$,  we extend $\phi$ to a new good partial coloring $\psi$ of $G$, in which $|A_{\psi}(e_{3})|\ge4$, $|A_{\psi}(e_{i})|\ge 5$ for $i=1,2,4,5$ and  $|A_{\psi}(f_{j})|\ge 3$ for $j=2,4$. Observe  that  $A_{\psi}(e_{5})\cap A_{\psi}(f_{2}) \neq \emptyset$, we can further  extend $\psi$ to another good partial coloring $\sigma$ of $G$ by coloring $e_{5}$ and $f_{2}$ with the same color $\beta\in A_{\psi}(e_{5})\cap A_{\psi}(f_{2})$. Then we have $|A_{\sigma}(e_{3})|\ge3$, $|A_{\sigma}(e_{i})|\ge 4$ for $i=1,2,4$ and  $|A_{\sigma}(f_{4})|\ge 2$.

Now, if $A_{\sigma}(e_{2})\cap A_{\sigma}(f_{4}) \neq \emptyset$, we can color $e_{2}$ and $f_{4}$ with the same color $\gamma\in A_{\sigma}(e_{2})\cap A_{\sigma}(f_{4})$ and then color $e_{3},e_{4},e_{1}$ by SDR. Otherwise, we greedily color $f_4, e_3,e_4,e_1,e_2$ in this order. 

{\bf Case 2}. $|E(\{u_1,u_1'\},\{u_2,u_2'\})|=0$.

In this case, $w_{1},w_{1}', w_{2}, w_{2}'$ are not necessarily distinct. However, this does not affect the following arguments. Please see Figure \ref{fig:degree2_case2} for the names of vertices and edges. Set $S=\{v_{0},v_{1},v_{2}\}$. By Lemma  \ref{greedy}, there is a good partial coloring $\phi$ of $G$ with six uncolored edges $e_1,e_2,f_1,f_2,f_3,f_4$. It is clear that $|A_{\phi}(e_{1})|=|A_{\phi}(e_{2})|=4$ and $|A_{\phi}(f_{j})|\ge 2$ for $j=1,2,3,4$.
  
  	\begin{figure}[htbp]
  	\centering
  	\resizebox{11.1cm}{3.5cm}{\includegraphics{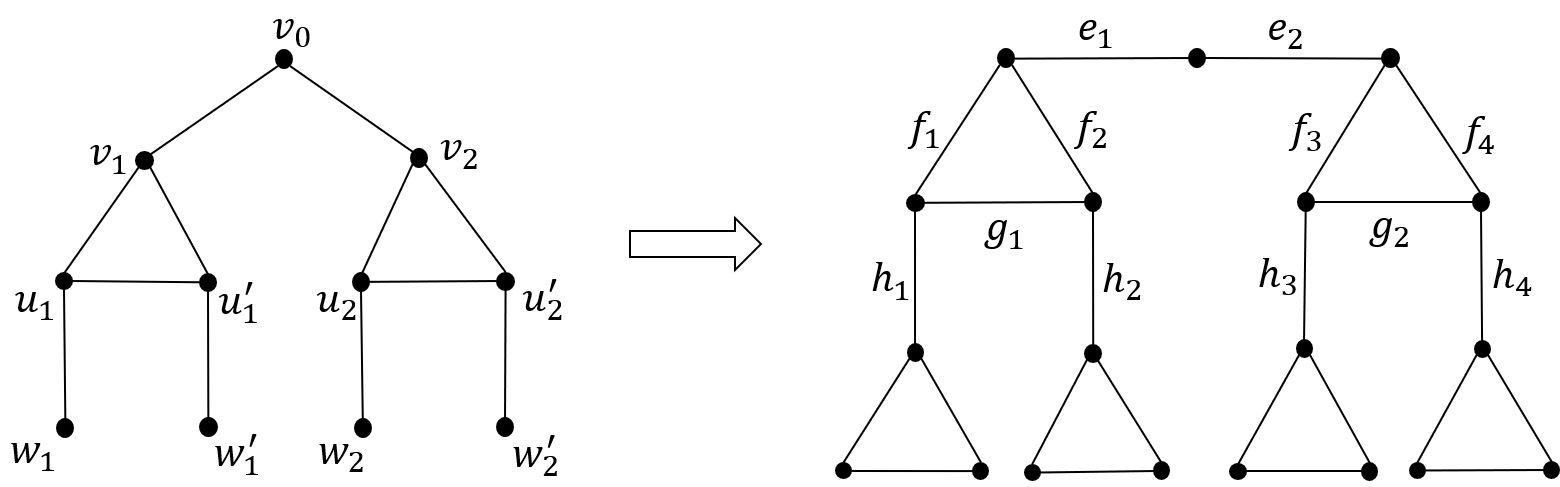}}
  	\caption{Case 2 in the proof of Lemma \ref{lemma:degree-2}}
  	\label{fig:degree2_case2}
       \end{figure}
  
We now make two easy but useful observations. One is that $A_{\phi}(f_{1})\cup A_{\phi}(f_{2})\subseteq A_{\phi}(e_{1})=[1,7]\setminus F_{\phi}(e_{1})$ and $A_{\phi}(f_{3})\cup A_{\phi}(f_{4})\subseteq A_{\phi}(e_{2})=[1,7]\setminus F_{\phi}(e_{2})$, where $F_{\phi}(e_{1})=\{\phi(g_{1}),\phi(h_{1}),\phi(h_{2})\}$ and $F_{\phi}(e_{2})=\{\phi(g_{2}),\phi(h_{3}),\phi(h_{4})\}$. The other is that, for any $i\in\{1,2\}$ and any $j\in\{3,4\}$, $f_{i}$ and $f_{j}$ do not see each other. 

As $|A_{\phi}(e_{1})|=|A_{\phi}(e_{2})|=4$, $A_{\phi}(e_{1})\cap A_{\phi}(e_{2})\neq \emptyset$. According to the value of $|A_{\phi}(e_{1})\cap A_{\phi}(e_{2})|$, we divide the proof into the following four subcases.

{\bf Subcase 2.1}.  $|A_{\phi}(e_{1})\cap A_{\phi}(e_{2})|=1$.  ($|A_{\phi}(e_{1})\cup A_{\phi}(e_{2})|=7$.)

In this case, let's explain why $\phi$ can be extend to a good coloring of $G$ by greedily coloring $f_{1},f_{2},f_{3},f_{4},e_{1}$ and $e_{2}$ in this order. First, $f_{1},f_{2},f_{3},f_{4}$  can be colored properly. Secondly, since  $A_{\phi}(f_{3})\cup A_{\phi}(f_{4})\subseteq A_{\phi}(e_{2})$ and $|A_{\phi}(e_{1})\cap A_{\phi}(e_{2})|=1$, there  is at least one color available  for $e_{1}$ when $e_{1}$ is to be colored. Finally, since $A_{\phi}(f_{1})\cup A_{\phi}(f_{2})\subseteq
A_{\phi}(e_{1})$  and  $|A_{\phi}(e_{1})\cap A_{\phi}(e_{2})|=1$, after $e_{1}$ is colored, at least one color is  available  for $e_{2}$.

 {\bf Subcase 2.2}. $|A_{\phi}(e_{1})\cap A_{\phi}(e_{2})|=2$. ($|A_{\phi}(e_{1})\cup A_{\phi}(e_{2})|=6$.)
 
Let $A_{\phi}(e_{1})\cap A_{\phi}(e_{2})=\{\alpha_{1},\alpha_{2}\}$. If $\{\alpha_{1},\alpha_{2}\}\cap (A_{\phi}(f_{1})\cup A_{\phi}(f_{2}))\cap (A_{\phi}(f_{3})\cup A_{\phi}(f_{4}))\neq\emptyset$, without loss of generality, assume $\alpha_{1} \in A_{\phi}(f_{1})\cap A_{\phi}(f_{3})$. Now, a good coloring of $G$ can be obtained by first coloring $f_{1}$ and $f_{3}$ with $\alpha_{1}$ and then greedily coloring $f_{2},f_{4},e_{1},e_{2}$ in this order. Thus, we assume that $\{\alpha_{1},\alpha_{2}\}\cap (A_{\phi}(f_{1})\cup A_{\phi}(f_{2}))\cap (A_{\phi}(f_{3})\cup A_{\phi}(f_{4}))=\emptyset$.

If $\{\alpha_{1},\alpha_{2}\}\cap (A_{\phi}(f_{1})\cup A_{\phi}(f_{2})) =\{\alpha_{1},\alpha_{2}\}\cap (A_{\phi}(f_{3})\cup A_{\phi}(f_{4})) =\emptyset$, then it is easy to see that $\phi$ can be extended to a good coloring of $G$. Therefore,  by symmetry, we may assume that  $\alpha_{1}\in A_{\phi}(f_{1})$ and  $\alpha_{1}\notin A_{\phi}(f_{3})\cup A_{\phi}(f_{4})$.

Suppose that  $A_{\phi}(f_{1})\cup A_{\phi}(f_{2})\neq \{\alpha_{1},\alpha_{2}\}$. If $\alpha_{1}\in A_{\phi}(f_{2})$, then there exists a color $\beta\in (A_{\phi}(f_{1})\cup A_{\phi}(f_{2})) \setminus \{\alpha_{1},\alpha_{2}\}$. We can  color  $f_{1}$ and $f_{2}$ to get a new good partial coloring $\psi$ of $G$ so that $\{\psi(f_{1}),\psi(f_{2})\}=\{\alpha_{1},\beta\}$. And then  we extend $\psi$  to  a good coloring of $G$ by  greedily coloring  $f_{3},f_{4},e_{2},e_{1}$ in this order. If $\alpha_{1}\notin A_{\phi}(f_{2})$, then there exists a color $\beta\in  A_{\phi}(f_{2}) \setminus \{\alpha_{1},\alpha_{2}\}$. A good coloring of $G$ can be obtained by coloring $f_{1}$ with $\alpha_{1}$,  $f_{2}$ with $\beta$ and then coloring the remaining edges in the order $f_{3},f_{4},e_{2},e_{1}$. 

Now we assume $A_{\phi}(f_{1})\cup A_{\phi}(f_{2}) = \{\alpha_{1},\alpha_{2}\}$. 
Recall that $A_{\phi}(f_{3})\cup A_{\phi}(f_{4}) \subseteq A_{\phi}(e_{2})$,
we must have $A_{\phi}(f_{3})= A_{\phi}(f_{4}) = A_{\phi}(e_{2})\setminus \{\alpha_{1},\alpha_{2}\}$.
Notice that $g_{1}$ sees a color $\alpha$ in $\phi$ if and only if $f_{1}$ or $f_{2}$ also sees $\alpha$ in $\phi$, we can recolor $g_{1}$ with $\alpha_{1}$ to obtain a new good partial coloring of $G$, which we  refer to it as $\psi$.  Observe that $A_{\psi}(f_{1})= A_{\psi}(f_{2}) = \{\phi(g_{1}),\alpha_{2}\}$, $A_{\psi}(e_{1})=(A_{\phi}(e_{1})\setminus \{\alpha_{1}\}) \cup  \{\phi(g_{1})\}$, $A_{\psi}(f_{3})= A_{\psi}(f_{4}) =A_{\phi}(f_{3})$  and  $A_{\psi}(e_{2})= A_{\phi}(e_{2})$.
Now, we can get a good coloring of $G$ by coloring $f_{1}$ with $\phi(g_{1})$, $f_{2}$ with $\alpha_{2}$,  $f_{3}$ and $f_{4}$ with the two colors in $A_{\phi}(e_{2})\setminus \{\alpha_{1},\alpha_{2}\}$, $e_{2}$ with $\alpha_{1}$ and  $e_{1}$ with a color in $A_{\phi}(e_{1})\setminus \{\alpha_{1},\alpha_{2}\}$.
 
{\bf Subcase 2.3}. $|A_{\phi}(e_{1})\cap A_{\phi}(e_{2})|=3$. ($|A_{\phi}(e_{1})\cup A_{\phi}(e_{2})|=5$.)

Let $A_{\phi}(e_{1})\cap A_{\phi}(e_{2})=\{\alpha_{1},\alpha_{2},\alpha_{3}\}$, $A_{\phi}(e_{1})=\{\alpha_{1},\alpha_{2},\alpha_{3},\beta\}$ and $A_{\phi}(e_{2})=\{\alpha_{1},\alpha_{2},\alpha_{3},\gamma\}$.
Recall that $A_{\phi}(f_{1})\cup A_{\phi}(f_{2})\subseteq A_{\phi}(e_{1})$ and 
$A_{\phi}(f_{3})\cup A_{\phi}(f_{4})\subseteq A_{\phi}(e_{2})$, as $|A_{\phi}(f_{j})|\ge 2$, we must have
$A_{\phi}(f_{j})\cap\{\alpha_{1},\alpha_{2},\alpha_{3}\}\neq \emptyset$ for each $j=1,2,3,4$.
We have two subcases to deal with.

{\bf Subcase 2.3.1}. $(A_{\phi}(f_{1})\cup A_{\phi}(f_{2}))\cap (A_{\phi}(f_{3})\cup A_{\phi}(f_{4}))\neq \emptyset$.

W.l.o.g., let $\alpha_{1} \in A_{\phi}(f_{1}) \cap A_{\phi}(f_{3})$.
By assigning the color $\alpha_{1}$ to $f_{1}$ and $f_{3}$ and coloring $f_{2}$ (resp. $f_{4}$) with a color in $A_{\phi}(f_{2})\setminus \{\alpha_1\}$ (resp. $A_{\phi}(f_{4})\setminus \{\alpha_1\}$), we extend $\phi$ to a new  good partial coloring $\psi$. Observe that $|A_{\psi}(e_{1})|\ge 1$, $|A_{\psi}(e_{2})|\ge 1$
and $|A_{\psi}(e_{1})\cup A_{\psi}(e_{2})|=|(A_{\phi}(e_{1})\cup A_{\phi}(e_{2}))\setminus\{\alpha_{1},\psi(f_{2}),\psi(f_{4})\}|\ge 2$. Therefore, we  can further extend $\psi$ to a good coloring of $G$ by SDR.

{\bf Subcase 2.3.2}. $(A_{\phi}(f_{1})\cup A_{\phi}(f_{2}))\cap (A_{\phi}(f_{3})\cup A_{\phi}(f_{4}))= \emptyset$.

In this subcase, it  is clear that $|A_{\phi}(f_{1})\cup A_{\phi}(f_{2})|=2$ or $|A_{\phi}(f_{3})\cup A_{\phi}(f_{4})|=2$.
W.l.o.g., we assume that $|A_{\phi}(f_{1})\cup A_{\phi}(f_{2})|=2$ (i.e., $A_{\phi}(f_{1})= A_{\phi}(f_{2})$), $\alpha_{1}\in A_{\phi}(f_{1})\cup A_{\phi}(f_{2})$ and $\alpha_{3}\in A_{\phi}(f_{3})\cup A_{\phi}(f_{4})$.

If $A_{\phi}(f_{1})= A_{\phi}(f_{2})= \{\alpha_{1},\alpha_{2}\}$, then $A_{\phi}(f_{3})= A_{\phi}(f_{4})= \{\alpha_{3},\gamma\}$. 
Recall that $g_{1}$ sees a color $\alpha$ in $\phi$ if and only if $f_{1}$ or $f_{2}$ also sees $\alpha$ in $\phi$.
Thus, we can always modify $\phi$ by recoloring $g_{1}$ with $\alpha_{1}$ to obtain a new good partial coloring $\psi$ of $G$, in which $A_{\psi}(f_{3})=A_{\psi}(f_{4})=\{\alpha_{3},\gamma\}$ and
$A_{\psi}(e_{2})=\{\alpha_{1},\alpha_{2},\alpha_{3},\gamma\}$.
Now, if $\phi(g_{1})=\gamma$,   we have 
$A_{\psi}(f_{1})=A_{\psi}(f_{2})=\{\alpha_{2},\gamma\}$ and
$A_{\psi}(e_{1})=\{\alpha_{2},\alpha_{3},\beta,\gamma\}$, and we are back to Subcase 2.3.1.
While if $\phi(g_{1})\neq\gamma$,  
we have 
$A_{\psi}(f_{1})=A_{\psi}(f_{2})=\{\alpha_{2},\phi(g_{1})\}$ and
$A_{\psi}(e_{1})=\{\alpha_{2},\alpha_{3},\beta,\phi(g_{1})\}$, implying $|A_{\psi}(e_{1})\cap A_{\psi}(e_{2})|=2$. And we  are back to Subcase 2.2.

If  $ A_{\phi}(f_{1})= A_{\phi}(f_{2})=\{\alpha_{1},\beta\}$,  by  symmetry,
we may  assume that $A_{\phi}(f_{3})\cup A_{\phi}(f_{4})\neq \{\alpha_{2},\alpha_{3}\}$
(i.e., $\{\alpha_{3},\gamma\} \subseteq A_{\phi}(f_{3})\cup A_{\phi}(f_{4})$).
Also, we can  recolor $g_{1}$ with $\alpha_{1}$ to get a new good partial coloring $\psi$ of $G$, in which $\{\alpha_{3},\gamma\}\subseteq A_{\psi}(f_{3})\cup A_{\psi}(f_{4})$ and
$A_{\psi}(e_{2})=\{\alpha_{1},\alpha_{2},\alpha_{3},\gamma\}$.
Now,  if $\phi(g_{1})=\gamma$,
we have $A_{\psi}(f_{1})=A_{\psi}(f_{2})=\{\beta,\gamma\}$ and
$A_{\psi}(e_{1})=\{\alpha_{2},\alpha_{3},\gamma,\beta\}$, and we are back to Subcase 2.3.1.
And if $\phi(g_{1})\neq\gamma$,
 then we have 
 $A_{\psi}(f_{1})=A_{\psi}(f_{2})=\{\beta,\phi(g_{1})\}$ and
 $A_{\psi}(e_{1})=\{\alpha_{2},\alpha_{3},\beta,\phi(g_{1})\}$.  
 As $|A_{\psi}(e_{1})\cap A_{\psi}(e_{2})|=2$, we are back to Subcase 2.2.
 
{\bf Subcase 2.4}. $|A_{\phi}(e_{1})\cap A_{\phi}(e_{2})|=4$. ($|A_{\phi}(e_{1})\cup A_{\phi}(e_{2})|=4$.)

Let $A_{\phi}(e_{1})= A_{\phi}(e_{2})=\{\alpha_{1},\alpha_{2},\alpha_{3},\alpha_{4}\}$.
Erasing the colors of $g_{1}$ and $g_{2}$ in  $\phi$  yields another good partial coloring  $\psi$ of $G$.

If $|A_{\psi}(g_{1})|= |A_{\psi}(g_{2})|=1$, then it is straightforward to check that $A_{\phi}(f_{1})\cap A_{\phi}(f_{2})=\emptyset$, $A_{\phi}(f_{1})\cup A_{\phi}(f_{2})=A_{\phi}(e_{1})$, $A_{\phi}(f_{3})\cap A_{\phi}(f_{4})=\emptyset$ and $A_{\phi}(f_{3})\cup A_{\phi}(f_{4})=A_{\phi}(e_{2})$.
Now based on  $\phi$, we can always  color $f_{1},f_{2},f_{3},f_{4}$  to obtain a  new good partial coloring $\sigma$ with $\{\sigma(f_{1}),\sigma(f_{2})\}=\{\sigma(f_{3}),\sigma(f_{4})\}$. Consequently, $e_{1}$ and $e_{2}$  can be colored properly.

Now, w.l.o.g.,  we assume that $|A_{\psi}(g_{1})| \ge 2$. In this situation, it is not difficult to verify that 
$A_{\psi}(g_{1})\cap \{\alpha_{1},\alpha_{2},\alpha_{3},\alpha_{4}\}\neq \emptyset$.  
W.l.o.g., let $\alpha_{1}\in A_{\psi}(g_{1})$.
By recoloring the edge $g_1$ with the color $\alpha_1$ in $\phi$, we obtain a new good partial coloring $\sigma$ of  $G$, in which $A_{\sigma}(e_{1})=\{\phi(g_{1}),\alpha_{2},\alpha_{3},\alpha_{4}\}$ and $A_{\sigma}(e_{2})=\{\alpha_{1},\alpha_{2},\alpha_{3},\alpha_{4}\}$.
Obviously, $|A_{\sigma}(e_{1}) \cap A_{\sigma}(e_{2})|=3$ and so we are back to Subcase 2.3.

As we have exhausted all cases, the lemma follows. 
\end{pf}

\begin{lemma}\label{lemma:cut vertex}
Let $G$ be a connected claw-free cubic graph. If $G$ has a cut vertex, then $\chi'_{s}(G)\leq 7$.
\end{lemma}

\begin{pf}	
	Let $v_0$ be a cut vertex  of  $G$ and $N(v_{0})=\{u_{0},v_{1},v_{2}\}$. 
	As $G$ is a claw-free cubic graph, there is exactly one edge in $G$ among the vertices in $N(v_{0})$.
	Without loss of generality, let $v_{1}v_{2}\in E(G)$. 
	Let $G_{1}$  be the component of $G-v_{0}$ containing $u_{0}$ and  $G_{2}$ the component of $G-v_{0}$ containing $v_{1}v_{2}$. Please see Figure \ref{fig:cut vertex} for the names of vertices and edges in $G$.
	
	\begin{figure}[htbp]
	\centering
	\resizebox{7.3cm}{5cm}{\includegraphics{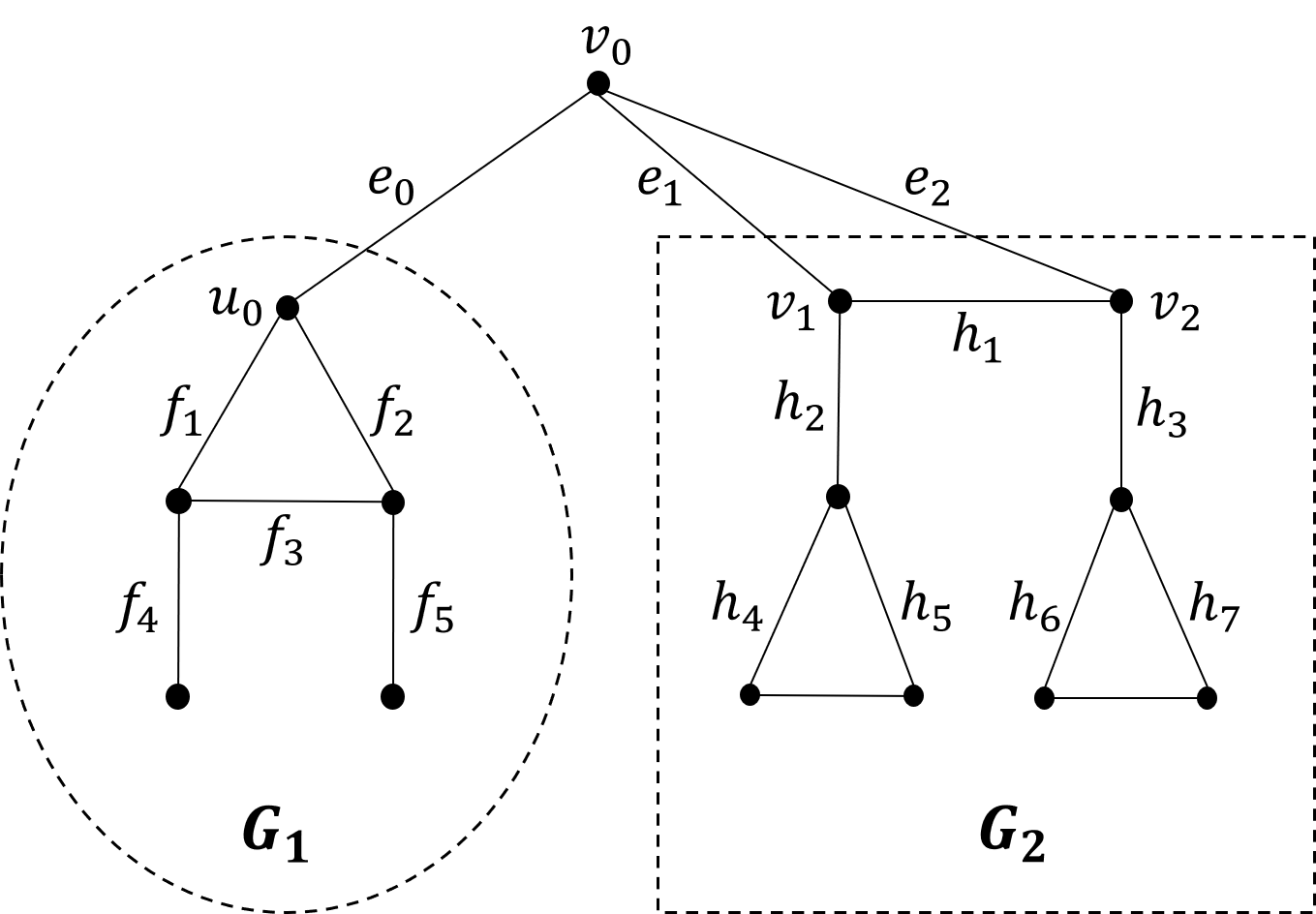}}
	\caption{A graph $G$ with a cut vertex $v_{0}$}
	\label{fig:cut vertex}
       \end{figure}
	
	By  Lemma \ref{lemma:degree-2}, $G_{1}$ (resp. $G_{2}$) has a good coloring, say $\phi_{1}$ (resp. $\phi_{2}$). 
	We may assume that $\{\phi_{1}(f_{1}),\phi_{1}(f_{2}),\phi_{1}(f_{3})\}=\{\phi_{2}(h_{1}),\phi_{2}(h_{2}),\phi_{2}(h_{3})\}$ as otherwise we can permute the colors among the edges  of $G_{1}$. 
Combining $\phi_{1}$ and $\phi_{2}$ yields a good partial coloring of $G$, calling it $\phi$.
	It is clear that $|A_{\phi}(e_{0})|=2$, $|A_{\phi}(e_{1})|\ge 2$ and $|A_{\phi}(e_{2})|\ge 2$.

If $|A_{\phi}(e_{0})\cup A_{\phi}(e_{1})\cup A_{\phi}(e_{2})|\ge 3$, then we can obtain a good coloring of $G$ by SDR.
Thus we assume that
 $A_{\phi}(e_{0})= A_{\phi}(e_{1})= A_{\phi}(e_{2})=\{\alpha_{1},\alpha_{2}\}$.  Let $\beta=\phi(f_{4})$. 
Permute the two colors $\alpha_{1}$ and $\beta$ in $E(G_{1})$.
 After the permutation,  we have $A_{\phi}(e_{0})=\{\alpha_{2},\beta\}$ and $A_{\phi}(e_{1})= A_{\phi}(e_{2})=\{\alpha_{1},\alpha_{2}\}$. 
 Since $|A_{\phi}(e_{0})\cup A_{\phi}(e_{1})\cup A_{\phi}(e_{2})|\ge 3$, we complete the proof  by SDR.
\end{pf}

\begin{lemma}\label{lemma:$4$-cycle with a chord}
	Let $G$ be a  $2$-connected claw-free cubic graph.
	If $G$ contains a $4$-cycle with a chord, then $\chi'_{s}(G)\le 7$.
\end{lemma}

\begin{pf}
	Let $C=v_{1}v_{2}v_{3}v_{4}v_{1}$ be a $4$-cycle in $G$ with $v_{1}v_{3}\in E(G)$.
	Because $G$ is not isomorphic to $K_{4}$,  $v_{2}v_{4}\notin E(G)$.
Let  $u_{2}$ and $u_{4}$ denote the third neighbor of $v_{2}$ and $v_{4}$, respectively.
Since $G$ is a $2$-connected claw-free cubic graph, it is not difficult to see that  $u_{2}\neq u_{4}$ and $u_{2}u_{4}\notin E(G)$.
If $N(u_{2})\cap N(u_{4})\neq \emptyset$, then $G$ is isomorphic to $H_{4}$ (as shown in Figure \ref{fig:H_{4}}) and so $\chi'_{s}(G)\le 7$.
	Thus we assume that $N(u_{2})\cap N(u_{4})= \emptyset$. Please see Figure \ref{fig:H_{4}} for the the names of vertices and edges in $G$. 
	
		\begin{figure}[htbp]
		\centering
		\resizebox{15cm}{4.3cm}{\includegraphics{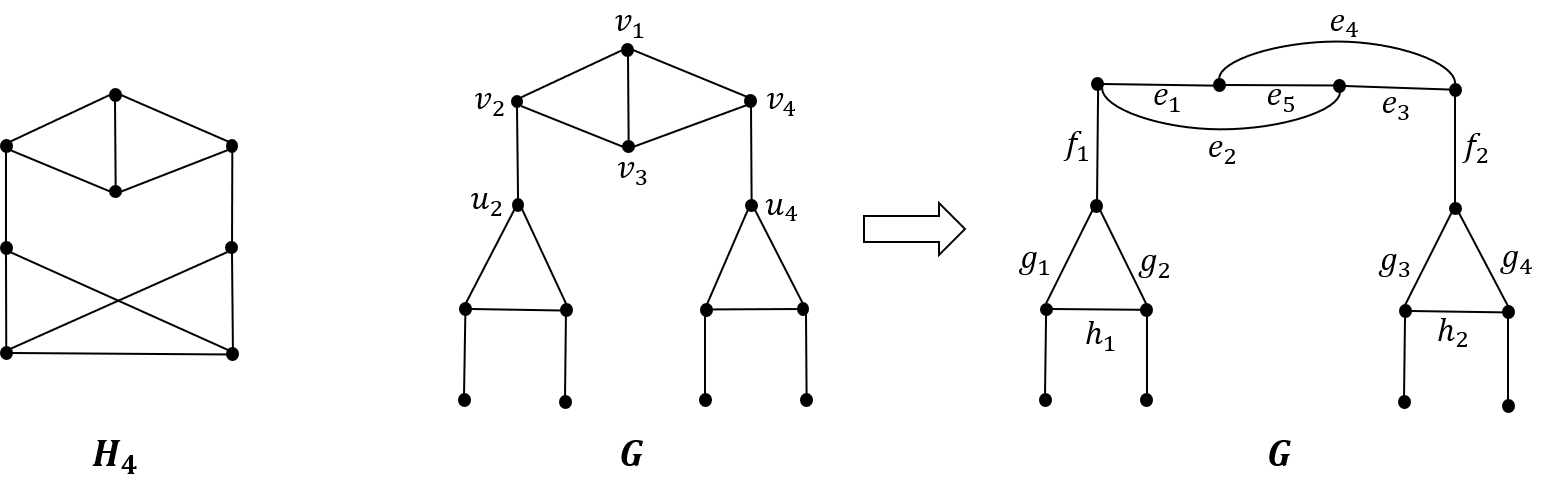}}
		\caption{The graph $H_{4}$ and a graph $G$ with a chorded $4$-cycle}
		\label{fig:H_{4}}
	       \end{figure}
	
	Set $S=\{v_{1},v_{2},v_{3},v_{4}\}$. By  Lemma  \ref{greedy}, we get a good partial coloring $\phi$ of $G$, in which
	$|A_{\phi}(e_{i})|=5$ for $i\in [1,4]$, $|A_{\phi}(e_{5})|=7$ and $|A_{\phi}(f_{1})|=|A_{\phi}(f_{2})|=2$.
	If $|A_{\phi}(e_{1})\cup A_{\phi}(e_{3})|\ge 6$, then a good coloring of $G$ can be easily obtained by 	SDR.
	Thus we assume that $|A_{\phi}(e_{1})\cup A_{\phi}(e_{3})|=5$. 
	Observe that $A_{\phi}(e_{1})=A_{\phi}(e_{2})$ and $A_{\phi}(e_{3})=A_{\phi}(e_{4})$, 
	we may assume that $A_{\phi}(e_{1})=A_{\phi}(e_{2})=A_{\phi}(e_{3}) =A_{\phi}(e_{4})=[1,5]$.
	And so $\{\phi(g_{1}),\phi(g_{2})\}=\{\phi(g_{3}),\phi(g_{4})\}=\{6,7\}$.
		
	Erasing the colors of  $g_{1}$ and $g_{2}$ in $\phi$ yields a good partial coloring of $G$, calling it  $\psi$. 
	Now, if  there exists some color $\alpha\in A_{\psi}(g_{1})\setminus \{6,7\}$,
	we  color $g_{1}$  with $\alpha$ and $g_{2}$ with $\phi(g_{2})$ to get a new good partial coloring $\sigma$ of $G$, in which $A_{\sigma}(e_{1})=([1,5]\setminus \{\alpha\})\cup \{\phi(g_{1})\}$ and $ A_{\phi}(e_{3})=[1,5]$. 
	Clearly, $|A_{\sigma}(e_{1})\cup A_{\sigma}(e_{3})|\ge 6$.
	So we can extend $\sigma$  to a good coloring of $G$ by SDR.
	Thus we assume $A_{\psi}(g_{1})= \{6,7\}$.
	Symmetrically, we may also assume $A_{\psi}(g_{2})=\{6,7\}$.
	Notice that $h_{1}$ sees a color $\alpha$ in $\psi$ if and only if $g_{1}$ or $g_{2}$ also sees  $\alpha$ in $\psi$,
	we can recolor $h_{1}$ with  $\phi(g_{1})$ in $\psi$. By further coloring $g_{1}$ with $\phi(h_{1})$ and  $g_{2}$ with $\phi(g_{2})$, we get a new good partial coloring $\sigma$ of $G$ with $A_{\sigma}(e_{1})=([1,5]\setminus \{\phi(h_{1})\})\cup \{\phi(g_{1})\}$ and $ A_{\phi}(e_{3})=[1,5]$. Again we have $|A_{\sigma}(e_{1})\cup A_{\sigma}(e_{3})|\ge 6$ and so  complete the proof by SDR.
\end{pf}

\begin{lemma}\label{lemma:induced 4-cycle}	
Let $G$ be a  $2$-connected  claw-free cubic graph.	
	If $G$ contains an induced  $4$-cycle, then $\chi'_{s}(G)\le 7$.
\end{lemma}

\begin{pf}
	Let $C=v_{1}v_{2}v_{3}v_{4}v_{1}$ be an induced $4$-cycle in $G$.
	As $G$ is a claw-free cubic graph not isomorphic to the $3$-prism, $u_{1}u_{2}\notin E(G)$.
	 Let $N(u_{1})=\{v_{1},v_{2},w_{1}\}$ and  
	$N(u_{2})=\{v_{3},v_{4},w_{2}\}$ (Please refer to Figure \ref{fig:induced 4cycle}).
	Since $G$ is  $2$-connected,   $w_{1}w_{2}\notin E(G)$.
	
	\begin{figure}[htbp]
		\centering
		\resizebox{12cm}{3.4cm}{\includegraphics{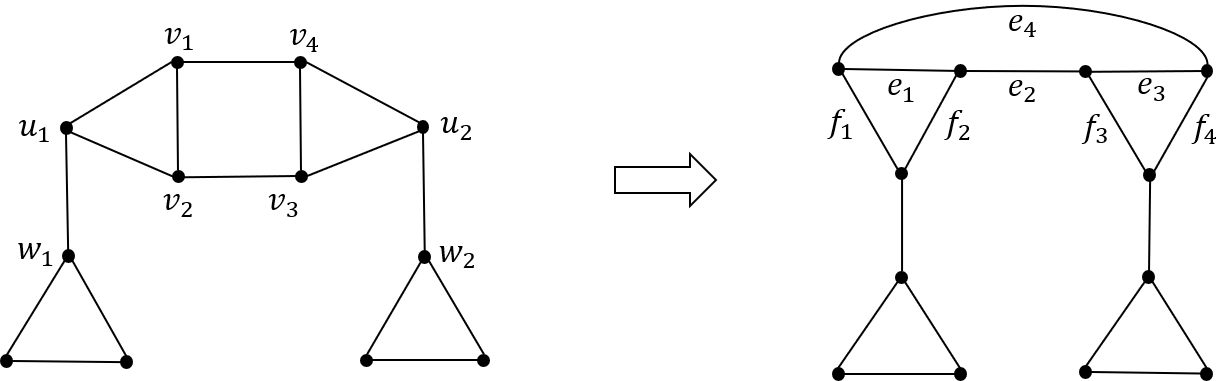}}
		\caption{A graph $G$ with an induced $4$-cycle $C=v_{1}v_{2}v_{3}v_{4}v_{1}$}
		\label{fig:induced 4cycle}
	\end{figure}

	Setting $S=\{v_{1},v_{2},v_{3},v_{4}\}$, by Lemma \ref{greedy}, we get a good partial coloring $\phi$ of $G$,  
	in which $|A_{\phi}(e_{i})|=6$ for $i=1,3$, 
       $|A_{\phi}(e_{i})|\ge 5$ for $i=2,4$ and
       $|A_{\phi}(f_{j})|=4$  for $j=1,2,3,4$.
Observe  that there exists some color $\alpha \in A_{\phi}(f_{1})\cap A_{\phi}(f_{3})$, 
by coloring $f_{1}$ and $f_{3}$ with $\alpha$, we extend $\phi$ to another good partial coloring $\psi$ of $G$, 
in which $|A_{\psi}(f_{2})|=|A_{\psi}(f_{4})|=3$,
 $|A_{\psi}(e_{1})|=|A_{\psi}(e_{3})|=5$ and 
$|A_{\psi}(e_{2})|=|A_{\psi}(e_{4})|\ge 4$.

Now, if there exists some $\beta \in A_{\psi}(f_{2})\cap A_{\psi}(f_{4})$,
 then we can color  $f_{2}$ and $f_{4}$ with $\beta$  and  greedily color $e_{2},e_{4},e_{1},e_{3}$ in this order  to 
get a good coloring of $G$.
Otherwise, we have 
$|A_{\psi}(f_{2})\cup A_{\psi}(f_{4})|=6$, 
and so the remainning  six egdes can be colored  properly by SDR.  Thus  the lemma holds.
\end{pf}

At present time, the only case we need to deal with is that  $G$ is a $2$-connected claw-free cubic graph without $4$-cycles. It is clear that each vertex is exactly on one $3$-cycle. Let $\Delta_{1}$ and  $\Delta_{2}$ be two $3$-cycles in $G$. Then $\Delta_{1}$ and  $\Delta_{2}$ are vertex-disjoint and there is at most one edge from $V(\Delta_{1})$  to  $V(\Delta_{2})$.  Moreover, it is not difficult to see that $G$ has no induced odd cycles of length at least $5$ and each induced even cycle has length at least $6$. 	

\begin{lemma}\label{lemma:induced 2p-cycle}
If $G$ is a $2$-connected claw-free cubic graph without $4$-cycles, then $\chi'_{s}(G)\le 7$.
\end{lemma}

 \begin{pf}
Choose a minimum induced even cycle  $C$ in $G$ and label the vertices of $C$ as 
$v_{1},v_{2},\ldots,v_{2p}$ and the edges $e_{i}=v_{i-1}v_{i}$ for $i=2,3,\ldots,2p$ and $e_{1}=v_{2p}v_{1}$.
 Clearly, $2p\ge6$.
 As $G$ is a claw-free cubic graph  and has no $4$-cycles, we may assume that $v_{2j-1}$ and $v_{2j}$ have the common neighbor $u_{2j}$ for each $j\in[1,p]$ and $u_{2},u_{4},\ldots,u_{2p}$ are distinct.
 Since $C$ is a minimum induced even cycle  in $G$,  $u_{2i}u_{2j}\notin E(G)$ for any $i,j\in[1,p]$.
 Let $w_{2j}$ denote the third neighbor of $u_{2j}$ for each $j\in[1,p]$.
 Also, $w_{2},w_{4},\ldots,w_{2p}$ are distinct because $G$ is claw-free and cubic.
 As for the names of vertices and edges in  $G$, please refer to Figure \ref{fig:induced 6cycle}.
 
 \begin{figure}[htbp]
 	\centering
 	\resizebox{16cm}{4.5cm}{\includegraphics{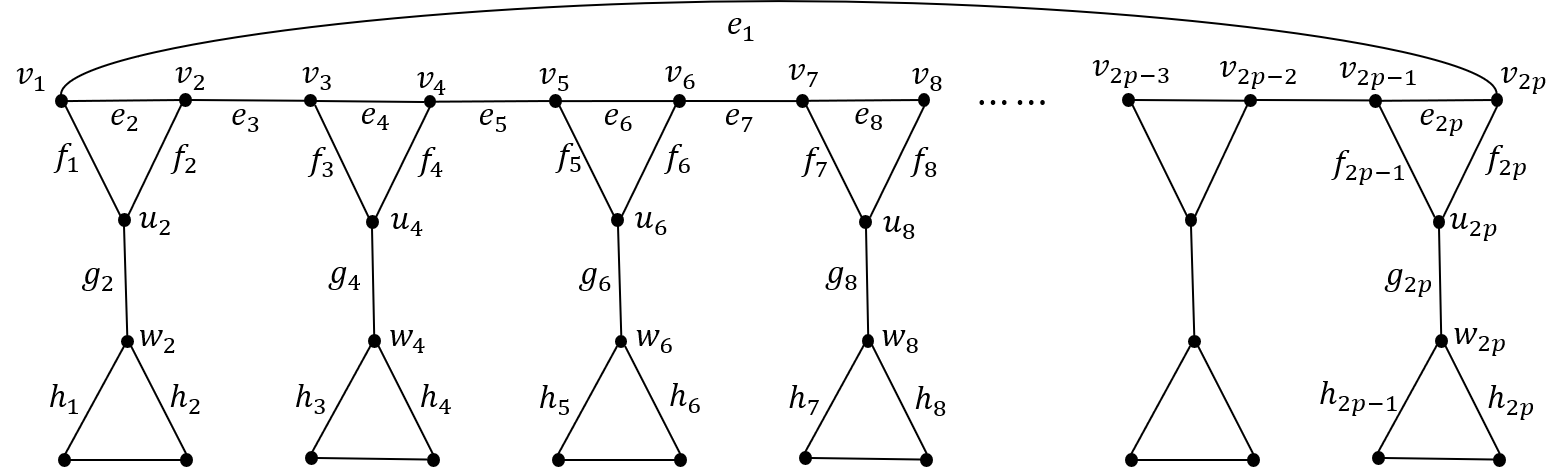}}
 	\caption{A graph $G$ with an induced $2p$-cycle $(p\ge3)$}
 	\label{fig:induced 6cycle}
 \end{figure}

Notice that if $2p\ge 8$, then $w_{2j}w_{2j+2}\notin E(G)$ for any $j\in[1,p-1]$ and $w_{2}w_{2p}\notin E(G)$.
And if $2p=6$, then there is at most one edge in the subgraph of $G$ induced by $\{w_{2},w_{4},w_{6}\}$  since $G$ is not isomorphic to $K_{4}^{\Delta}$.
So in either cases, we may assume that  $w_{2}w_{4}\notin E(G)$. This implies that the four edges $h_{1},h_{2},h_{3},h_{4}$ are distinct.

Let $G'$ be the graph obtained from $G$ by deleting $v_{1},v_{2},\ldots,v_{2p},u_{2},u_{4},u_{6}$ 
and adding a new edge $w_{2}w_{4}$.
Observe that $G'$  is a claw-free subcubic graph and each component of $G'$ has a vertex of degree less than $3$,  by Lemmas \ref{lemma:degree-1} and \ref{lemma:degree-2}, $G'$  has a good coloring $\phi$.  
Ignoring $w_{2}w_{4}$ in $\phi$ yields a good partial coloring of $G$.
Let $\alpha=\phi(w_{2}w_{4})$. And let $c_{i}=\phi(h_{i})$ for each $i\in[1,2p]$.
Then $\alpha,c_{1},c_{2},c_{3},c_{4}$ are distinct. 
As $|A_{\phi}(g_{6})|\ge2$,  there exists a color $\beta\in A_{\phi}(g_{6})\setminus\{\alpha\}$.
Without loss of generality, assume that  $c_{3}\neq \beta$.

Now based on $\phi$, we color $g_{2},g_{4}$ and $e_{6}$ with the same color $\alpha$, 
color $g_{6}$ with $\beta$ and $e_{5}$ with $c_{3}$. 
Observe that $f_{5}$ sees at most $5$ colors  at this moment,  
we can color $f_{5}$ with some color $\gamma\neq c_{4}$. 
This yields a new  coloring $\psi$, which is indeed a good partial coloring of $G$. 
Please refer to Figure \ref{fig:Tn1} and  Figure \ref{fig:Tn2} for this coloring.

 \begin{figure}[htbp]
	\centering
	\resizebox{16cm}{4.5cm}{\includegraphics{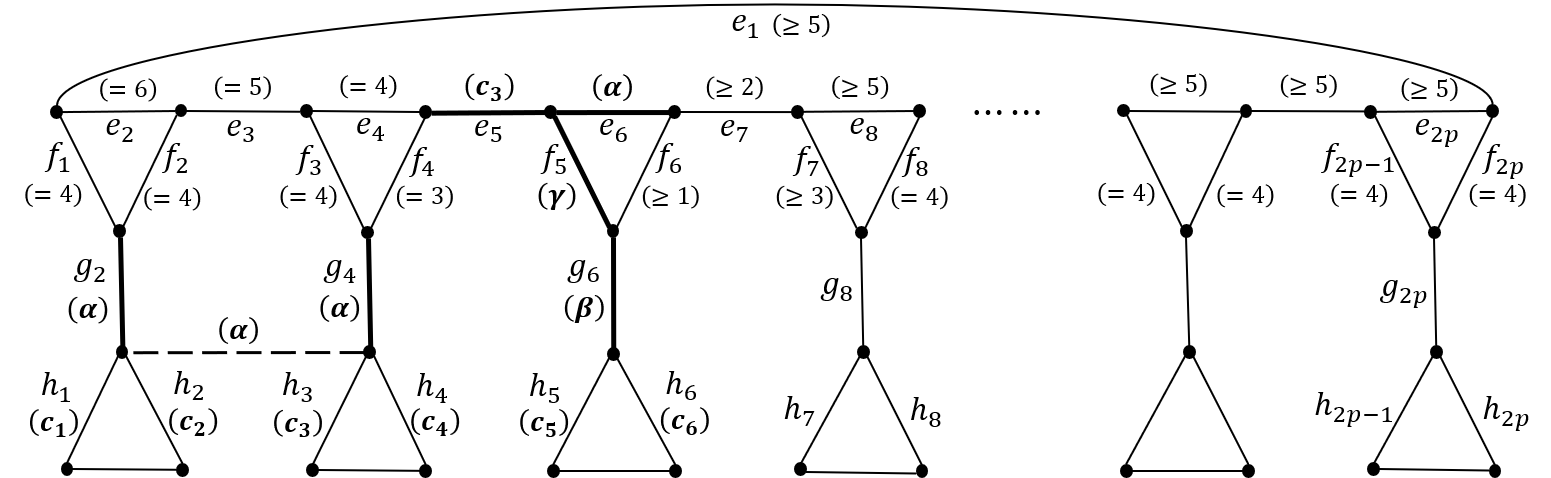}}
	\caption{$\psi$, $\bar{G}_{\psi}$ and 
		$| A_{\psi}(e)|$ for each $e\in\bar{E}_{\psi}$ when $2p\ge 8$}
	\label{fig:Tn1}
 \end{figure}
 
  \begin{figure}[htbp]
	\centering
	\resizebox{8cm}{4.5cm}{\includegraphics{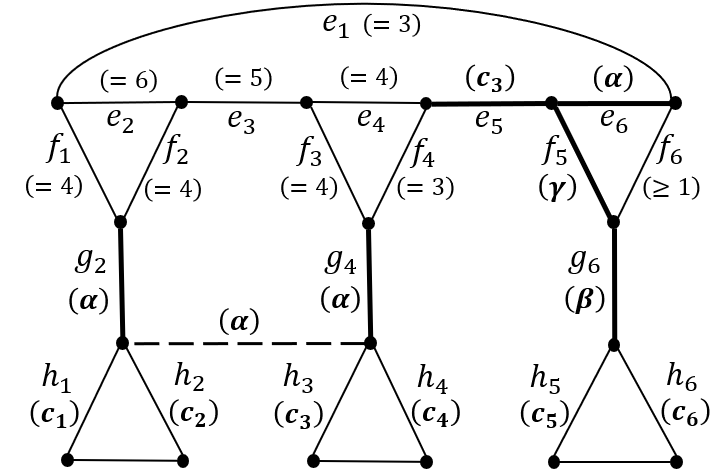}}
	\caption{$\psi$, $\bar{G}_{\psi}$ and 
		$| A_{\psi}(e)|$ for each $e\in\bar{E}_{\psi}$ when $2p=6$}
	\label{fig:Tn2}
 \end{figure}

In the following,  we will construct  a good coloring $\sigma$ of  $\bar{G}_{\psi}$
so that $\sigma(e)\in A_{\psi}(e)$ for every $e\in\bar{E}_{\psi}$. 
Then, this coloring $\sigma$ combining with the coloring $\psi$ forms a good coloring of $G$, and we are done.

If $2p\ge 8$, then it is straightforward to check that 
$|A_{\psi}(e_{1})|\ge 5$, $|A_{\psi}(e_{2})|=6$, $|A_{\psi}(e_{3})|=5$, 
$|A_{\psi}(e_{4})|=4$, $|A_{\psi}(e_{7})|\ge 2$, $|A_{\psi}(e_{i})|\ge 5$ for each $i\in \{8,9,\dots,2p\}$, 
$|A_{\psi}(f_{4})|=3$, $|A_{\psi}(f_{6})|\ge 1$, $|A_{\psi}(f_{7})|\ge 3$, 
and $|A_{\psi}(f_{i})|=4$ for each $i\in \{1,2,3,8,9,\dots,2p\}$. (Please refer to Figure \ref{fig:Tn1}). 
And if $2p=6$,  it is clear that 
$|A_{\psi}(e_{1})|=3$, $|A_{\psi}(e_{2})|=6$, $|A_{\psi}(e_{3})|=5$, $|A_{\psi}(e_{4})|=4$,
$|A_{\psi}(f_{4})|=3$, $|A_{\psi}(f_{6})|\ge 1$, and $|A_{\psi}(f_{i})|=4$ for each $i\in \{1,2,3\}$.
 (Please refer to Figure \ref{fig:Tn2}).

Before we go further to the next step, we make two observations which are simple but helpful in the following arguments (please refer to Figure \ref  {fig:Property} for the illustrations of Observation A and Observation B).\\
 
\noindent{\bf Observation A}:  $A_{\psi}(f_1)=A_{\psi}(f_2)=[1,7]\setminus \{\alpha,c_1,c_2\}\subseteq A_{\psi}(e_2) = [1,7]\setminus \{\alpha\}$.\\

\noindent{\bf Observation B}:  $A_{\psi}(e_4)=A_{\psi}(f_4)\cup \{c_4\}$, $A_{\psi}(f_3)=A_{\psi}(f_4)\cup \{\gamma\}$, and $A_{\psi}(e_3)=A_{\psi}(f_4)\cup \{c_4,\gamma\}$. 
 \begin{figure}[htbp]
	\centering
	\resizebox{15.5cm}{2.8cm}{\includegraphics{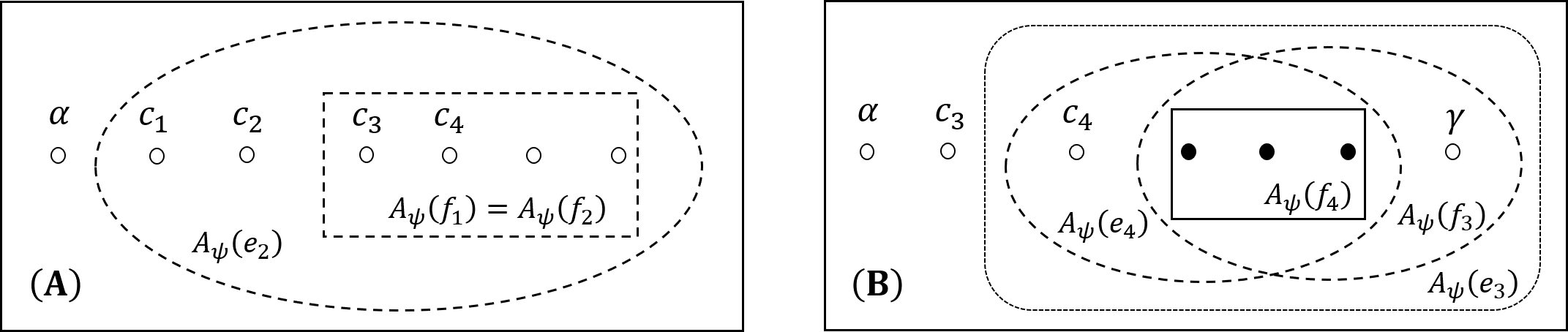}}
	\caption{Observetion A and Observetion B}
	\label{fig:Property}
\end{figure} 
 
Now, if $2p\ge 8$, we color  $f_{6},e_{7},f_{7},e_{8},f_{8} \ldots,e_{2p},f_{2p},e_{1},f_{1},e_{2},f_{2},e_{3},f_{3},f_{4}$ in this order greedily 
to get a good partial  coloring $\sigma$ of  $\bar{G}_{\psi}$, where only $e_{4}$ is uncolored.
And if $2p=6$, we color $f_{6},e_{1},f_{1},e_{2},f_{2},e_{3},f_{3},f_{4}$  in this order greedily 
to get a good partial coloring $\sigma$ of  $\bar{G}_{\psi}$, with only $e_{4}$ being uncolored.
For simplicity, let $a_{i}=\sigma(e_{i})$ and  $b_{j}=\sigma(f_{j})$ for edges $e_{i},f_{j}\in \bar{E}_{\psi}$ that are colored in $\sigma$ (please refer to Figure \ref{fig:Tn1_color0}).

 \begin{figure}[htbp]
	\centering
	\resizebox{16cm}{4.5cm}{\includegraphics{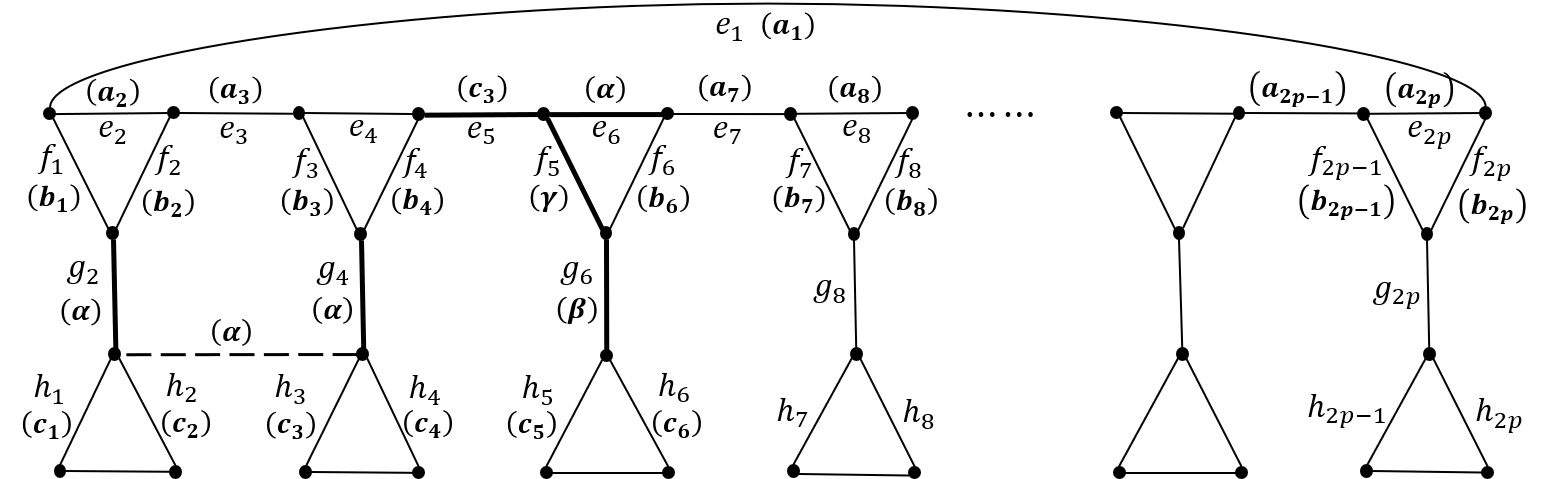}}
	\caption{the good partial coloring $\sigma$ of $\bar{G}_{\psi}$}
	\label{fig:Tn1_color0}
 \end{figure}

If $A_{\psi}(e_4)\setminus \{a_2,b_2,a_3,b_3,b_4\}\not=\emptyset$, then the edge $e_4$ can be colored properly and we are done. We therefore assume that  $A_{\psi}(e_4)\subseteq \{a_2,b_2,a_3,b_3,b_4\}$.

If $\gamma\not\in \{a_2,b_2,a_3\}$, then we may assume that $b_3=\gamma$ as otherwise we can  recolor the edge $f_3$ with the color $\gamma$  to make $b_3=\gamma$. By Observation B, $\gamma\not\in A_{\psi}(f_4)$ and $\gamma\not\in A_{\psi}(e_4)$. It follows that there exists a color $b_4^*$ in $A_{\psi}(f_4)\setminus \{a_3,b_3,b_4\}$ and $A_{\psi}(e_4)=\{a_2,b_2,a_3,b_4\}$. Now, we can recolor $f_4$ with $b_4^*$ and then color $e_4$ with $b_4$. This results in a good coloring of $\bar{G}_{\psi}$. Therefore we assume that $\gamma\in \{a_2,b_2,a_3\}$.  

It follows from the assumption $\gamma\in \{a_2,b_2,a_3\}$ that $b_3\not=\gamma$. 
Since  $b_3\in A_{\psi}(f_3)$, by Observation B, $b_3\in A_{\psi}(f_4)$. 
Therefore we must have $\{b_3,b_4\}\subset A_{\psi}(f_4)$.
The remainder of the proof is divided into two cases according to 
whether $a_3$ is in $A_{\psi}(f_4)$ or not.

{\bf Case 1.} $a_3\not\in A_{\psi}(f_4)$.

Since $a_3\in A_{\psi}(e_3)$ and, by Observation B, $A_{\psi}(e_3)=A_{\psi}(f_4)\cup \{c_4,\gamma\}$, it is clear that $a_3\in \{c_4,\gamma\}$. Let $A_{\psi}(f_4)=\{b_3,b_4, \alpha^*\}$. Then $A_{\psi}(e_4)=\{c_4,b_3,b_4, \alpha^*\}$. If  $c_4$ (or $\alpha^*$) is not contained in $\{a_2,b_2,a_3\}$, then we can color the edge $e_4$ with $c_4$ (or $\alpha^*$) and obtain a good coloring of $\bar{G}_{\psi}$. 
Thus we assume that $\{c_4,\alpha^*\} \subset \{a_2,b_2,a_3\}$. This together with $\gamma\in \{a_2,b_2,a_3\}$ implies that  $\{a_2,b_2,a_3\}=\{c_4,\alpha^*,\gamma\}$. As $a_3\in\{c_4,\gamma\}$, $\alpha^*$ must be in $\{a_2,b_2\}$. Now we can recolor the edge $f_4$ with $\alpha^*$ and then color $e_4$ with $b_4$. This completes the coloring of $\bar{G}_{\psi}$.

{\bf Case 2.} $a_3\in A_{\psi}(f_4)$. (So $A_{\psi}(f_4)=\{a_3,b_3,b_4\}$.)

In this case, by Observation B, $A_{\psi}(e_4)=A_{\psi}(f_4)\cup \{c_4\}=\{c_4,a_3,b_3,b_4\}$ and $a_3\not=\gamma$. Since $A_{\psi}(e_4)\subseteq \{a_2,b_2,a_3,b_3,b_4\}$, $c_4$ must be in $\{a_2,b_2\}$.
 As $\gamma \in \{a_2,b_2,a_3\}$ and $a_3\not=\gamma$, $\gamma\in\{a_2,b_2\}$.
Therefore, $\{a_2,b_2\}=\{c_4,\gamma\}$.
Now we erase the colors of $f_{3}$ and  $f_{4}$ in $\sigma$  to get a new good  partial coloring of $\bar{G}_{\psi}$, where the only uncolored edges are $f_{3}$, $f_{4}$ and $e_{4}$. 
This  coloring is still denoted by $\sigma$.
We next extend $\sigma$ to a good coloring of $\bar{G}_{\psi}$.

If $a_{2}\in A_{\psi}(f_{2})$, then, by Observation A, we also have  $a_{2}\in A_{\psi}(f_{1})$.
Notice that $b_{1},a_{2}\notin \{a_{2p},b_{2p},a_{1},b_{2},a_{3}\}$,  we exchange the colors of $f_{1}$ and $e_{2}$ in $ \sigma$.
If $a_{2}=c_{4}$  and $b_{2}=\gamma$ (so $b_{1}\neq c_{4}$),  we  color $e_{4}$  with $c_{4}$. 
It follows that there is at least one color available for $f_{3}$ and 
there are at least two colors available for $f_{4}$.
Hence we can color $f_{3}$ and  $f_{4}$ properly. 
If  $a_{2}=\gamma$  and $b_{2}=c_{4}$  (so $b_{1}\neq \gamma$), we color $f_{3}$ with  $\gamma$. 
Then $e_{4}$ has at least one available color  and  $f_{4}$ has at least two available colors.
Therefore, coloring $e_{4}$ and then $f_{4}$ gives a good coloring of $\bar{G}_{\psi}$. 

If $a_{2}\notin A_{\psi}(f_{2})$,  
then there exists some color $b_{2}^{*}$ in $A_{\psi}(f_{2})\setminus \{a_{1},b_{1},b_{2}\}$. 
We  recolor $f_{2}$ with $b_{2}^{*}$.
Now, if $a_{2}=c_{4}$ and $b_{2}=\gamma$, we color $f_{3}$ with $\gamma$.  
Recall that $\gamma\not\in A_{\psi}(e_4)$, 
there is at least one color available for $e_{4}$ and there are at least two colors available for $f_{4}$.
Thus $e_4$ and $f_4$ can be colored properly.
If $a_{2}=\gamma$ and $b_{2}=c_{4}$,  we  recolor $e_{3}$ with $c_{4}$.
Then  both  $e_{4}$ and $f_{3}$ have at least two available colors  and  $f_{4}$ has three available colors.
It follows that we can obtain a  good coloring of $\bar{G}_{\psi}$ by SDR. 
 \end{pf}
 
Combining all lemmas in this section, Theorem \ref{main} holds.

\section{Final Remarks}

Let $G$ be any connected claw-free subcubic graph not isomorphic to the triangular prism.
This paper proves that $\chi'_{s}(G)\leq 7$.  And this upper bound $7$ is sharp.
In addition, our proof implies a linear-time
algorithm for finding a strong $7$-edge-coloring for such a graph.

If $G$ is a connected claw-free cubic graph not isomorphic to the triangular prism, then it is easy to see that $\chi'_{s}(G)\geq 6$. Therefore, for such graph $G$, $\chi'_{s}(G)\in\{6,7\}$.  Let $H$ be a connected cubic graph. Denote by $H^{\Delta}$ the graph obtained from $H$ by replacing each vertex with a $3$-cycle. It is clear that $H^{\Delta}$ is a connected claw-free cubic graph. 

We end this paper by asking the following three questions.\\

\noindent {\bf Question 1}: Let $k\ge 3$ be an integer and $H$ the $k$-prism. Is it true that $\chi'_{s}(H^{\Delta})=6$?\\

\noindent {\bf Question 2}: Is it possible to characterize all connected claw-free cubic graphs $G$ with $\chi'_{s}(G)=6$?\\

\noindent {\bf Question 3}: Is the  strong list-chromatic index of any claw-free subcubic graph other than the triangular prism at most 7?\\

\noindent{\bf Declaration of competing interest}

The authors declare that they have no known competing financial interests or personal relationships that could have appeared to influence the work reported in this paper.




\begin{thebibliography}{99}
	
\bibitem{A1992} L.D. Andersen,   The strong chromatic index of a cubic graph is at most 10,  Discrete Math. 108  (1992) 231-252.

\bibitem{BBH2015} J.  Bensmail,  M. Bonamy,  H. Hocquard,
Strong edge coloring sparse graphs,  Electron. Notes Discret. Math. 49 (2015) 773-778.

\bibitem{BPP2022} M. Bonamy,  T. Perrett,  L. Postle,
Colouring graphs with sparse neighbourhoods: Bounds and applications, J. Comb. Theory, Ser. B 155 (2022) 278-317.

\bibitem{BJ2015} H. Bruhn, F. Joos, A stronger bound for the strong chromatic index,  Electron. Notes Discret. Math. 49 (2015) 277-284.

\bibitem{C2006} D. Cranston, Strong edge-coloring of graphs with maximum degree 4 using 22 colors, Discrete Math. 306 (21) (2006) 2772-2778.

\bibitem{DJS2020} M. D{e}bski,  K. Junosza-Szaniawski,  M. {\'S}leszy{\'n}ska-Nowak,
Strong chromatic index of $K_{1,t}$-free graphs, Discrete Appl. Math.  284 (2020) 53-60. 

\bibitem{E1988} P. Erd\H{o}s, Problems and results in combinatorial analysis and graph theory, Discrete Math. 72 (1988) 81–92.

\bibitem{EN1989} P. Erd\H{o}s,  J. Ne\v{s}et\v{r}il,
Irregularities of Partitions, Springer,  (1989) 162-163.

\bibitem{FSGT1990} R.J. Faudree,  R.H. Schelp,   A. Gy{\'a}rf{\'a}s,  Zs. Tuza,  The strong chromatic index of graphs, Ars Comb. 29 (1990) 205-211.

\bibitem{FJ1983} J.L. Fouquet,  J.L. Jolivet, 
Strong edge-colorings of graphs and applications to multi-$k$-gons, Arc Comb. 16 (1983) 141-150.

\bibitem{H1935} P. Hall, On representatives of subsets, J. Lond. Math. Soc. 10 (1935) 26-30.

\bibitem{H1990}  P. Hor{\'a}k, The strong chromatic index of graphs with maximum degree four,
 Contemp. Methods Graph Theory 399 (1990) 403.
 
\bibitem{HQT1993} P. Hor{\'a}k,  H. Qing,  W.T.  Trotter, Induced matchings in cubic graphs,  J. Graph Theory 17 (2) (1993) 151-160.
 
\bibitem{HSY2018} M. Huang,  M. Santana, G. Yu,
Strong chromatic index of graphs with maximum degree four, Electron. J. Comb. 25 (3) (2018) \#P3.31.

\bibitem{HdK2021} E. Hurley,  R. de Joannis de Verclos,  R.J. Kang, An improved procedure for colouring graphs of bounded local density, In Proceedings of the 2021 ACM-SIAM Symposium on
Discrete Algorithms (SODA), pages 135–148, 2021. https://arxiv.org/abs/2007.07874.

\bibitem{KLRSWY2016} A.V. Kostochka,  X. Li,  W. Ruksasakchai,  M. Santana,  T. Wang,  G. Yu,  Strong chromatic index of subcubic planar multigraphs, 
Eur. J. Comb. 51 (2016) 380-397.

\bibitem{LLY2018} J.-B. Lv,  X. Li,  G. Yu, On strong edge-coloring of graphs with maximum degree 4, Discrete Appl. Math. 235 (2018) 142–153. 

\bibitem{LLZ2022}  J.-B. Lv,  J. Li,  X. Zhang, 
On strong edge-coloring of claw-free subcubic graphs, 
Graphs Comb.  38 (3) (2022) 63.

\bibitem{MR1997} M. Molloy,  B. Reed, A bound on the strong chromatic index of a graph, J. Comb. Theory, Ser. B  69 (2) (1997) 103-109.

\bibitem{NKGB2000} T. Nandagopal,  T. Kim, X. Gao,   V. Bharghavan, Achieving MAC layer fairness in wireless packet networks, Proceedings of the 6th annual international conference on Mobile computing and networking, (2022) 87-98.

\bibitem{R1997} S. Ramanathan, A unified framework and algorithm for (T/F/C) DMA channel assignment in wireless networks, Proceedings of INFOCOM'97, IEEE, 
2 (1997) 900-907.

\bibitem{WSWC2018} Y. Wang, W.C. Shiu, W. Wang,  M. Chen, Planar graphs with maximum degree 4 are strongly 19-edge-colorable, Discrete Math. 341 (6) (2018) 1629-1635.

\bibitem{WL2008} J. Wu, W. Lin,   The strong chromatic index of a class of graphs, Discrete Math.  308 (2008) 6254-6261.

\bibitem{Z2015}  C. Zang, The strong chromatic index of graphs with maximum degree $\Delta$, (2015), 
arXiv preprint arXiv:1510.00785.

\end{thebibliography}
\end{document}